\documentclass[11pt]{article}
\usepackage{localeng}
\usepackage[margin=25mm]{geometry}                

\newtheorem*{theorem}{Theorem}
\newtheorem*{lemma}{Lemma}

\title{Ergodic theorem and algorithmic randomness}
\author{Alexander Shen\thanks{\raggedright LIRMM, University of Montpellier, CNRS, Montpellier, France, \texttt{alexander.shen@lirmm.fr}, \texttt{sasha.shen@gmail.com}. Supported by FLITTLA
ANR-21-CE48-0023 grant. Author is grateful to all members of the Kolmogorov seminar, especially to Vladimir Vyugin for patient explanations of his argument.}}
\date{}

\begin{document}
\maketitle

\begin{abstract}
We prove the constructive version of Birkhoff's ergodic theorem following Vyugin~\cite{vyugin} but trying to separate and state explicitly the combinatorial statement on which this proof is based. We pose some questions related to this statement (and the effective ergodic theorem in general).

\end{abstract}

\section{Birkhoff's theorem}

Let $T$ be a measure-preserving mapping of a probability space $X$ into itself. Let $f\colon X\to \mathbb{R}$ be a bounded integrable function. Consider the average of $f$ along $T$-trajectories, i.e., the limit
     $$
 \lim_{n\to\infty}\frac{f(x)+f(Tx)+f(TTx)+\ldots + f(T^{(n-1)}(x))}{n}    
     $$
(that may exist or not). Birkhoff's theorem claims that \emph{this limit exists for almost all $x$}. In addition, \emph{if the transformation $T$ is ergodic} (i.e., every set $X'\subset X$ that almost everywhere coincides with  $T^{-1}(X')$, has measure $0$ or $1$), then this limit is (almost everywhere) equal to the average of $f(x)$ over all points~$x$.

Note that the second statement is an easy corollary of the first one. Indeed, the value of this limit (considered as a function of $x$) is $T$-invariant. Therefore for every constant $c$ the set of all points where the limit is smaller than $c$ is either a set of measure zero or a set of full measure. When the parameter $c$ increases, we should jump from the first case to the second one, so there is some threshold $c_0$ where this jump happens. We conclude that the limit equals $c_0$ almost everywhere, and use Lebegue's convergence theorem for bounded functions to conclude that $c_0$ is the average of $f$ over~$X$.

Bishop~\cite{bishop} suggested a constructive version of Birkhoff's theorem that can be proved in the framework of his constructivism program. Using his ideas (in particular the upcrossing inequalities), Vyugin~\cite{vyugin} proved that the statement of Birkhoff's theorem is valid not just for \emph{almost all $x$}, but for \emph{all algorithmically random $x$}. In our note we provide an exposition of the (classical) Birkhoff's theorem using Bishop -- Vyugin arguments that states explicitly the core (pure combinatorial) lemma on which this proof is based. Then we explain (assuming some acquaintance with algorithmic randomness) why this argument proves Vyugin's result.

\section{How Birkhoff's theorem can be proven}

Recall the statement of the theorem: the set of all $x$ such that the sequence
     $$
f(x),\  \frac{f(x)+f(Tx)}{2},\ \frac{f(x)+f(Tx)+f(TTx)}{3},\ldots     \eqno(*)
     $$
has no limit, is a null set (has measure~$0$). This sequence is bounded since the function~$f$ is bounded. If the sequence has no limit, then it crosses back and forth some interval $(\alpha,beta)$ with $\alpha<\beta$ infinitely many times. We may consider only rational intervals, then there are countably many of them. So it is enough to prove that for every ``gap'' $(\alpha,\beta)$ \emph{the set of all sequences such that $(*)$ crosses $(\alpha,\beta)$ infinitely many times has measure $0$}. (By crossing we mean that some term of $(*)$ is smaller than $\alpha$, some later term is greater than $\beta$, then some term is again smaller than $\alpha$, etc.)
     
Let us fix some gap $(\alpha,\beta)$. We need to show that the measure of the set of point $x$ such that the sequence $(*)$ crosses this gap many times, is small. This would be done if we find an upper bound for the average number of crossings~$\mathbb{E}_x C(x)<\infty$. Here $C(x)$ is the number of crossing for $(*)$ (for given $x$) and the average is taken over~$x$. The function $C(x)$ is non-negative and can have infinite values; the finite upper bound means, in particular, that the set of $x$ where $C(x)$ is infinite has measure~$0$ (as we stated).

To get an upper bound for $\mathbb {E}_x C(x)$, we recall that $T$ is measure-preserving mapping. This allows us to add another layer of averaging \emph{over trajectories} (in addition to averaging over $x\in X$), and then change the order of averaging. In this way we see that an upper bound for average over trajectories is enough, and this bound could be proven for every bounded sequence. More precisely, we count the number of crossings for finite sequences (of arbitrarily large length) and prove an upper bound for the average that does not depend on the length. Let us explain now how it can be done.

Let us consider an arbitrary sequence $\mathbf{a}=(a_1, a_2, a_3,\ldots)$; we assume that $|a_i|\le A$ for some $A$ and all $i$. (Then we apply our conclusions to the sequence of $f$-values along a $T$-orbit.) Consider the sequence of averages
\[
E(\mathbf{a})=\left(a_1, \frac{a_1+a_2}{2}, \frac{a_1+a_2+a_3}{3}, \ldots, \frac{a_1+\ldots+a_k}{k},\ldots\right)
\]
Fix an arbitrary interval $(\alpha,\beta)$ and count how many times this sequence of averages crosses $(\alpha,\beta)$.

\begin{center}
\includegraphics[scale=1]{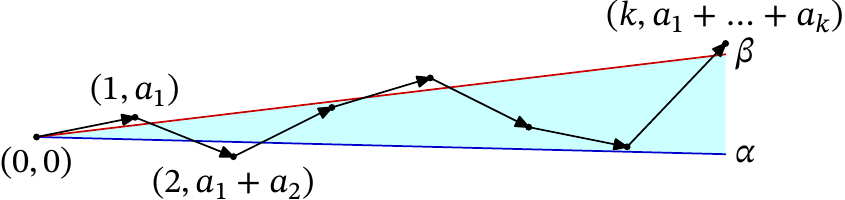}
\end{center}

It is instructive to draw a polygonal line starting from the origin; at each step we move right by $1$ and up/down by $a_i$. Then the average $(a_1+\ldots+a_i)/i$ is a slope of the line that goes from the origin to the end of the $i$-th segment. We are interested in the moments when this slope becomes less than $\alpha$ or greater than $\beta$, i.e., when the line crosses the rays with slopes $\alpha$ and $\beta$ (show in the picture). In our specific case (shown in the picture) we see one ``downcrossing'' (the second segment) and one ``upcrossing'' (third plus fourth segments). We have no more crossings in the picture (since the slope does not go below $\alpha$ after that).

The difference between upcrossing and downcrossing is at most $1$; let us count upcrossings. Formally, we consider increasing sequences of indices $l_1 < r_1 < l_2 < r_2<\ldots$ grouped into pairs $(l_i,r_i)$ such that $l_i$-th average is less than $\alpha$ and $r_I$-th average is greater than $\beta$ for all $i$. By $C_n(\mathbf{a})$ we denote the maximal number of pairs with these properties among the first $n$ terms of the sequence of averages $E(\mathbf{a})$. In our example we have $C_3(\mathbf{a})=0$ and$C_4(\mathbf{a})=\ldots=C_7(\mathbf{a})=1$. As $n$ grows, the number $C_n(\mathbf{a})$ (for a given $\mathbf{a}$) may only increase (when new upcrossing pairs appear). The limit value of this number (as $n\to infinite$), finite or infinite, is denoted by $C(\mathbf{a})$. We use the same letter $C$ as before, because our previous $C(x)$ is equal to $C(\mathbf{x})$ for $\mathbf{x}=(f(x), f(Tx), f(TTx),\ldots)$.

The following combinatorial lemme provides an upper bound for the number of upcrossing in a sliding window on a bounded sequence $\mathbf{a}=(a_1,a_2,\ldots)$. Let 
\[
\mathbf{a}^{(i)}=(a_i,a_{i+1}, a_{i+2},\ldots);
\]
compute the average number of upcrossing for $N$ neighbor positions of the window, i.e.,
\[
\frac{C_n(\mathbf{a}^{(1)})+C_n(\mathbf{a}^{(2)})+\ldots+C_n(\mathbf{a}^{(N)})}{N}.
\]

\begin{lemma}
For $N\ge n$ this average number does not exceed $\frac{A+|\alpha|+|\beta|}{\beta-\alpha}$ where $c$ is some absolute constant and $A$ is the upper bound for all $|a_i|$.
\end{lemma}

The condition $N\ge n$ guarantees that we have enough numbers for averaging. The denominator is $\beta-\alpha$; indeed, when the gap is narrow, we may have more upcrossings. The numerator contains $A$; when $A$ is large, the slope can be bigger and more upcrossings are possible. Note that multiplying $A$, $\alpha$ and $\beta$ by the same factor, we do not change the upper bound --- not surprising since multiplying all $a_i$, $\alpha$ and $\beta$ by the same constant, we do not change the number of upcrossings.

When shown in the picture, the switch from $\mathbf{a}^{(i)}$ to $\mathbf{a}^{(i+k)}$ means that the origin is moved $k$ steps to the right.
\begin{center}
\includegraphics[scale=1]{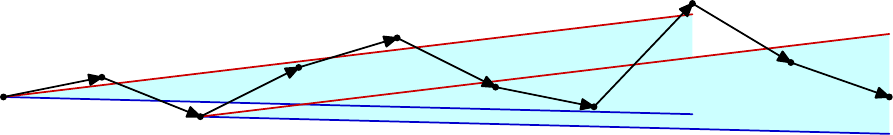}
\end{center}
This picture shows, for example, two triangles for which we measure the number of upcrossing: the initial one ($1$ upcrossing) and after a shift by $2$ (no upcrossings). Our lemma then provides an upper bound for the average number of upcrossings for $N$ triangles drawn along our polygonal line.

Using the Lemma, we now finish the proof or the ergodic theorem in the following way. We are interested in the average value of $\mathbb{E}_x C(x)$ where $C(x)=C(\mathbf{x})$ is the number of upcrossing for the sequence $\mathbf{x}=(f(x), f(Tx), f(TTx),\ldots)$. We need to show that this average is finite, and therefore the set of points $x$ where $C(x)$ is infinite has measure zero. To show this, it is enough to show that the average value of $C_n(\mathbf{x})$ (taken over $x$) is bounded by a constant that does not depend on $n$ (and then refer to the monotone convergence theorem). 

Here is the key observation that we already mentioned: \emph{since $T$ is measure-preserving, we have the same averages for $C_n(Tx)=C(\mathbf{x}^{(2)})$, for $C_n(TTx)=C(\mathbf{x}^{(3)})$ and so on}. So we add another layer of averaging over $N$ iterations of $T$. Then we change the order of averaging and note that now the internal average along the trajectory is bounded due to the lemma. (We may take $N$ greater than $n$, since the averaging argument is valid for every $N$.) So the average of  $C_n(\mathbf{x})$ over all $x$ does not exceed $c(F+|\alpha|+|\beta|)/(\beta-\alpha)$ where $F$ is the upper bound for the absolute value of $F$, and $(\alpha,\beta)$ is the gap interval we consider. 

This finishes the proof of the ergodic theorem (assuming the lemma is true; we prove it in the next section).

\section{Proof of the lemma}

We need to provide an upper bound for the sum of the numbers of upcrossings for several shifts of an arbitrary bounded sequence. We use the same idea that is often used for analysis of amortized algorithms. Namely, we introduce some \emph{potential function} and prove that it increases significantly when the number of upcrossings is large while remaining bounded.

This potential function (Vyugin used the name \emph{cumulative sum}) may look artificial, so let us start with several simple remarks. 

As Wikipedia says, ``In cycling, hiking, mountaineering and running, cumulative elevation gain refers to the sum of every gain in elevation throughout an entire trip''. For example, here
\begin{center}
\includegraphics[scale=1]{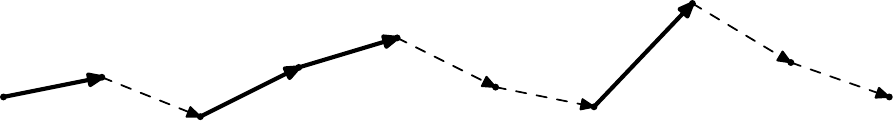}
\end{center}
we have three intervals when we go up, and the cumulative elevation gain is the sum of the height differences for these three intervals. 

Looking for a rigorous mathematical definition of cumulative elevation gain , we can use different approaches. For example, we may consider the minimal representation of our function (of bounded variation) as a difference of two non-decreasing functions, and consider the increase in the first one. However, for us the following definition will be useful:
\begin{center}
\includegraphics[scale=1]{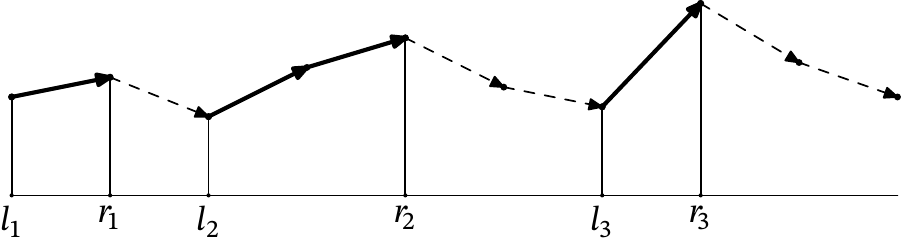}
\end{center}
we take all possible sequences of subintervals
\[
(l_1,r_1), (l_2,r_2),\ldots, (l_i,r_i)\qquad \text{with}\qquad 0\le l_1 < r_1 <  l_2 < r_2<\ldots<l_i<r_i\le n
\]
of our interval $[0,n]$; for each of them we consider the sum
\[
h(r_1)-h(l_1)+h(r_2)-h(l_2)+\ldots+h(r_i)-h(l_i),
\]
Then we take the maximum of this expression over all subsequences of intervals for all $i$ (including the case $i=0$, when this expression equals zero). This will indeed be the cumulative elevation gain (note that we may put all the increasing parts into disjoint intervals and avoid all decreasing parts).

We use this definition of the cumulative gain since it can be adapted for our purposes. Recall that we have fixed some gap $(\alpha,\beta)$. For all right endpoints $r_i$ we consider the height over the upper boundary of the gap, i.e., $h(r_i)-\beta r_i$, and for all left endpoints we consider $h(l_i)-\alpha l_i$. In other words, we consider the maximal value of the expression
\[
[h(r_1)-\beta r_1] - [h(l_1)-\alpha l_1]+[h(r_2)-\beta r_2] - [h(l_2)-\alpha l_2]+\ldots+[h(r_i)-\beta r_i]-[h(l_i)-\alpha l_i],
\]
The rest of the definition is unchanged (we still consider all families of disjoint intervals, and take the maximal value).

As we have announced, two things are of interest to us:
\begin{itemize}
\item the change of the cumulative sum when we shift the origin (in other words, delete the left segment in the polygonal line and add one more segment at the right end), and its connection to the number of upcrossings;
\item the lower and upper bounds for the cumulative sum.
\end{itemize}

Let us recall the notation. We start with a sequence $\mathbf{a}=(a_1,a_2,\ldots)$ and construct a polygonal line that goes from $(0,0)$ to $(1,a_1)$, then to $(2,a_1+a_2)$,\ldots, $(n, a_1+\ldots+a_n)$. We are interested in the number of upcrossings for the gap $(\alpha, \beta)$ and interval $[0,n]$ (left $n$ segments of the polygonal line), and also in the cumulative sums for intervals $[0,n]$ and $[1,n+1]$. More precisely, we compare the cumulative sums for the sequence $a_1,\ldots,a_n$ and $a_2,\ldots,a_{n+1}$ (where $a_{n+1}$ is the slope of the added segment).

What happens when we switch from $[0,n]$ to $[1,n+1]$? Let us look at the family of intervals $0\le l_1 < r_1 <l_2 <r_2<\ldots<l_s<r_s\le n$ that provides maximal value for the first cumulative sum. Let us assume for a while that $l_1>0$. Then all intervals $[l_i,r_i]$ can be used also for the second cumulative sum (with new origin, i.e., with endpoints decreased by~$1$). The corresponding terms become slightly bigger (since the quantity $\beta r_i$ that is subtracted at the right endpoint decreases by $\beta$ while $\alpha l_i$ decreases only by $\alpha$); each segment provides increase $(\beta -\alpha)$. So the total increase in the cumulative sum is at least $s(\beta-\alpha)$ where $s$ is the number of intervals in the maximal family. (The increase may be bigger if some other family of intervals will give a bigger value for $[1,n+1]$.)

Now consider the borderline case when $l_1=0$. Then the first interval $[l_1,r_1]$ should be decreased (the new interval is $[0,r_1-1]$) or should be abandoned completely (if $r_1=1$). But these changes can change the value of the cumulative sum only by $O(A+|\alpha|+|\beta|)$ where $A$ is the upper bound for all $|a_i|$; indeed, the changed zone of the graph has this size. 

So we conclude that the cumulative sum increases at least by
\[
s(\beta-\alpha)-O(A+|\alpha|+|\beta|),
\]
for the case when its maximal value is achieved for the family containing $s$ intervals.

Now let us show that a large number of upcrossings means that the cumulative sum expression reaches its maximal value on the large family of intervals (the number that we denoted by $s$ is large). Namely, we claim that the number of upcrossings is $O(s)$.  Recall that in the upcrossings count the points where the graph is below $\alpha$-line and the points where the graph is above $\beta$-line interleave. Consider, for example, three points $u<v<w$ where $u$ and $w$ are points of the first type, and $v$ is the point of the second type:
\begin{center}
\includegraphics[scale=1]{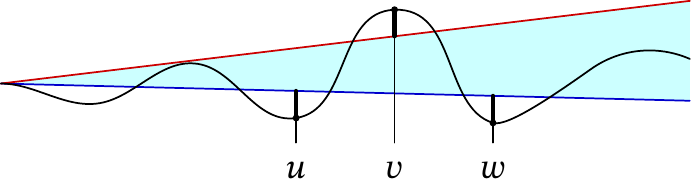}
\end{center}
(Here the line is shown as a curve.) Let us compare them with points $l_1<r_1<l_2<r_2<\ldots< l_s<r_s$ that provide the maximal value for the cumulative sum expression. The latter $2s$ points split the axis into intervals, and we claim that
\begin{center}
\emph{three points $u,v,w$ may not be inside the same interval}.
\end{center}
Indeed, imagine that they all are inside $(l_i,r_i)$; then this interval can be replaced by two intervals $[l_i,v]$ and $[w,r_i]$. This increases the expression for the cumulative sum, since we add two positive terms (that correspond to the vertical segments shown): one measures how much $v$-point is above the top line, and the other shows how much $w$ is below the bottom line. 

Similarly, if all three points $u,v,w$ are (strictly) inside $(r_i,l_{i+1})$, we can add new interval $[u,v]$ between $[l_i,r_i]$ and $[l_{i+1},r_{i+1}]$, and this also increases the expression for the cumulative sum (for the same reasons). 

Essentially the same argument works for the case when $u$- and $w$-points are above the top line and $v$-point is below the bottom line. Now we see (as we claimed) that the number of intervals in the maximal set of intervals is at least the number of upcrossing (up to some constant factor that can be easily computed, but for us any constant will work).

Now let us look at the bounds for the cumulative sum. The lower bound: as we mentioned, the cumulative sum is always non-negative (recall that we may consider an empty family of intervals). For the upper bound let us consider first the case $\alpha\ge 0$. Then the cumulative sum does not exceed the cumulative elevation gain. Indeed, for each interval $(i,j)$ we subtract $\beta j$ for the right endpoint and $\alpha i$ for the left endpoint, and $\beta j \ge \alpha j \ge \alpha i$. Therefore, for $\alpha\ge 0$ the cumulative sum is bounded by $nA$. On the other hand, the cumulative sum remains unchanged if we increase $\alpha$, $\beta$ and all terms $a_i$ by the same amount, so we may assume without loss of generality that $\alpha\ge 0$ after increasing $A$ by $|\alpha|$. Therefore, the cumulative sum for $n$ terms does not exceed  $(A+|\alpha|+|\beta|)n$  (here $\beta$ is added just for symmetry).

Now we may collect all the bounds together by adding the lower bounds for the increase of the cumulative sum for $N$ consecutive intervals of length $n$ with $c_1,\ldots,c_N$ upcrossings:
\[
\Omega(1)\sum_{i=1}^{N} c_i (\beta-\alpha)-NO(A+|\alpha|+|\beta|) \le O((A+|\alpha|+|\beta)n)
\]
Dividing by $N$ and recalling the assumption $n\le N$ we conclude that the average of $N$ upcrossing numbers $c_1,\ldots,c_N$ for $N$ consecutive intervals of length $n$ is 
 $O((A+|\alpha|+|\beta|)/(\beta-\alpha))$. This finishes the proof of the lemma.

\section{Effective versions}

It happens often that a classical theorem that says that something is true almost everywhere has an effective version saying that the same condition is true for all algorithmically random sequences. Usually some additional assumptions about effectivity of the objects appearing in the statement of the theorem are needed, and algorithmic randomness is understood in the sense of Martin-L\"of definition. This is also the case for the ergodic theorem. We restrict ourselves to the Cantor space of  infinite binary sequences and the shift operation (discarding the first term of a sequence). 

A measure $P$ on the Cantor space is determined by the measures of intervals $[x]$ (here $[x]$ , for a string $x$, is a set of all sequences that have prefix $x$), i.e., we may identify $P$ with the function $p(x)=P([x])$ defined on finite strings and taking non-negative relative values that satisfies two following conditions:
\begin{itemize}
\item $p(\Lambda)=1$;
\item $p(x)=p(x0)+p(x1)$;
\end{itemize} 
The measure $P$ is shift-invariant (stationary) if, in addition to these two conditions, the function $p$ satisfies one more condition:  
\begin{itemize}
\item $p(x)=p(0x)+p(1x)$.
\end{itemize}

Some effectivity assumptions are now needed. We assume that the measure $P$ is computable, i.e., there exists an algorithm that computes the values $p(x)$ with any requested precision. We also assume that the function $f$ is computable in the same way (the point $x$ is given to $f$ in the form of the oracle that provides approximations to $x$). The computability of $f$ implies its continuity; since the Cantor space is compact, every computable function is bounded.

\begin{theorem}[Vyugin~\cite{vyugin}]
 Let $P$ be an arbitrary computable stationary measure on the Cantor space. Let $f$ be a computable function on infinite sequences. Then for every Martin-L\"of random sequence $\omega$ there exists the average of $f$ along the orbit of $\omega$.
\end{theorem}

In other words, the set of sequences where the limit (average) does not exists is an effectively null set in the sense of Martin-L\"of: for every rational $\varepsilon>0$ we can effectively enumerate intervals that cover this set and have total measure at most $\varepsilon$.

We can take the characteristic function of an interval $[x]$ as $f$; then we conclude that for every random (with respect to $P$) sequence the limit frequency of factors $x$ exists. So zeros and ones appear in a $P$-random sequence with some limit frequencies, the same is true for two-bit blocks, etc.

How can we prove this effective version of the ergodic theorem? We may note that for every rational gap interval $(\alpha,\beta)$ the set of sequences that have more than $m$ upcrossings for the values of $f$ on their shifts, is effectively open. (Recall that $f$ can be computed with arbitrary precision, and we use strict inequalities ``less than $\alpha$'' and ``greater than $\beta$ in the definition of upcrossings. The measure of this effectively open set is bounded by $O(1/m)$ where the hidden constant depends only on the bound for $|f(\omega)|$ and the gap. 

In fact, this argument proves a stronger statement about sequences that are random with respect to the effectively closed class of stationary measures (the definition and basic properties can be found in~\cite{survey}):

\begin{theorem}[Vyugin, class randomness version]
Let $f$ be a computable function on the Cantor space. Then, for every sequence $\omega$ that is random with respect to the class of stationary measures, there exist an average of $f$ for the shifts of $\omega$.
\end{theorem}

This statement says that for every rational $\varepsilon>0$ we can enumerate a family of intervals that covers all bad sequence (=sequences for which the limit does not exist) that has total measure at most $\varepsilon$ for every stationary measure $P$. Note that this family of intervals \emph{does not depend on $P$} (and the measure $P$ may not be computable). This is exactly what is achieved by our construction (counting upcrossing we do not use any measure).

\section{Questions}

1. Is there a more direct proof of the lemma that does not use the cumulative sums (whose definition is quite artificial)?

One could start with the following special case. Let $X$ be some binary string. Consider the frequency of ones in its prefixes and its oscillations. For example, we may fix some thresholds, say $40\%$ and $60\%$ and count the number of oscillations of the frequency (how many times it goes below $40\%$ and then above $60\%$). Let us call this number ``oscillation number of $X$''. Then consider the following statement: for every $n$, if we take all factors of length $n$ in a sufficiently long string, the average oscillation number for these factors does not exceed $1000$ (the constant is chosen arbitrarily, but seems to be large enough). Can one prove this statement without cumulative sums? Note that the number of oscillations \emph{in one factor} can be arbitrarily large (if $n$ is large), this large number should be compensated by the other factors.

2. The statement of the ergodic theorem is valid for functions with values in $\mathbb{R}^2$ (and is a direct consequence of the statements for each coordinate). However, for two-dimensional sequences it is hard to define the upcrossing number in a natural way. Can we prove the two-dimensional results without considering two coordinates separately, in a more invariant way?

3. In fact, for the ergodic case (the shift is an ergodic transformation for the measure) the ergodic theorem can be extended to lower semicomputable functions $f$ (or for frequency of points in an effectively open set), see~\cite{hoyrup}. Can we combine these two results and prove convergence for non-ergodic case and lower semicomputable functions?

4. The ergodic theorem for amenable groups (in the special case of $\mathbb{Z}^2$) says that for a shift invariant measure $P$ on the set of all two-dimensional configurations of zeros and ones the set of all configurations that do not have the limit frequency of ones (over the increasing centered squares) has $P$-measure zero. Can we extend the effective version to this case and construct an (effectively?) open set that contains all those (limitless) sequences and has small measure with respect to any shift-invariant measure $P$? It seems that it is hard to adapt the upcrossing argument to this case.
\

\end{document}


\maketitle

\begin{abstract}
Мы излагаем конструктивное доказательство эргодической теоремы Биркгофа, предложенное в~\cite{vyugin}, выделив лежащую в его основе комбинаторную лемму. Формулируются некоторые вопросы, связанные с этим доказательством и этой леммой.
\end{abstract}

\section{Теорема Биркгофа}

Пусть $X$~--- вероятностное пространство, $T\colon X\to X$~--- измеримое отображение, сохраняющее меру, а $f\colon X\to \mathbb{R}$~--- ограниченная интегрируемая функция.  Тогда можно рассмотреть среднее функции $f$ по $T$-траектории, то есть предел
     $$
 \lim_{n\to\infty}\frac{f(x)+f(Tx)+f(TTx)+\ldots + f(T^{(n-1)}(x))}{n}    
     $$
Теорема Биркгофа утверждает, что \emph{этот предел существует при почти всех $x$}. Кроме того, \emph{если преобразование $T$ эргодично}, то есть не существует нетривиального (не нулевой и не полной меры) множества, которое (с точностью до меры $0$)  является собственным прообразом при преобразовании $T$, \emph{то этот предел \textup(почти всюду\textup) равен среднему функции $f$ по~$X$}.

Второе легко следует из первого, так как значение этого предела~--- инвариантная функция относительно преобразования $T$, и потому множество тех точек, где предел меньше $c$, может иметь лишь меру $0$ или $1$ (в силу эргодичности), и при пересечении параметром $c$ какой-то границы мера множества должна резко меняться с нуля на единицу. Это граничное значение равно среднему по теореме о мажорируемой сходимости (здесь используется ограниченность $f$).

Бишоп~\cite{bishop} предложил некоторый вариант теоремы Биркгофа, который можно доказать в рамках его конструктивного подхода к математике. Используя его идеи (upcrossing inequalities), Вьюгин~\cite{vyugin} доказал, что утверждение теоремы Биркгофа (при естественных предположениях) верно не просто для почти всех точек, но \emph{для всех алгоритмически случайных точек}.  В этой заметке мы приводим доказательство (классической) теоремы Биркгофа методом Бишопа -- Вьюгина, выделяя лежащую в его основе чисто комбинаторную лемму, а затем (предполагая знакомство с понятием алгоритмической случайности) указываем, как это рассуждение доказывает результат Вьюгина.

\section{Схема доказательства теоремы Биркгофа}

Теорема утверждает, что множество тех $x$, для которых последовательность
     $$
f(x),\  \frac{f(x)+f(Tx)}{2},\ \frac{f(x)+f(Tx)+f(TTx)}{3},\ldots     \eqno(*)
     $$
не имеет предела, имеет меру нуль. Эта последовательность ограничена (поскольку $f$ ограничена). Тогда отсутствие предела означает, что члены этой последовательности пересекают какой-то интервал $(\alpha,\beta)$ с $\alpha<\beta$ туда и обратно бесконечно много раз. Можно брать интервалы с рациональными концами, их счётное число, поэтому достаточно доказать, что для любой <<щели>> $(\alpha,\beta)$ \emph{множество тех $x$, для которых последовательность $(*)$ пересекает щель бесконечно много раз} (то есть бывает меньше $\alpha$, потом больше $\beta$, потом снова меньше $\alpha$ и так далее), \emph{имеет меру нуль}.

Зафиксируем щель $(\alpha,\beta)$. Как доказать, что мера множества тех $x$, для которых последовательность $(*)$ пересекает щель много раз, мала? Достаточно оценить среднее значение числа пересечений и показать, что $\mathbb{E}_x C(x)<\infty$. Здесь $C(x)$ --- число пересечений щели (если начать с~$x$), а математическое ожидание соответствует усреднению по~$X$. (Функция $C$ может быть для некоторых $x$ бесконечна, и наше утверждение о конечности математического ожидания означает, в частности, что множество таких $x$ имеет меру нуль --- как нам и хотелось.) 

Откуда возьмётся оценка для $\mathbb {E}_x C(x)$? Надо воспользоваться тем, что преобразование $T$ сохраняет меру: это позволит дополнить среднее \emph{по $X$} усреднением \emph{по траекториям}, а потом переставить порядок усреднения и доказать аналогичную оценку для \emph{средних по траекториям}.  А эта оценка уже будет верна для любой ограниченной последовательности. Чтобы сделать это аккуратно, надо сначала считать пересечения на конечном (но сколь угодно большом) участке и оценивать их среднее (величиной, не зависящей от длины участка). Вот как это делается.

Удобно отдельно рассмотреть колебания средних для произвольной последовательности (забыв об её происхождении как последовательности значений функции $f$ на орбите преобразования $T$). Пусть $\mathbf{a}=(a_1, a_2, a_3,\ldots)$ --- произвольная последовательность действительных чисел, не превосходящих по модулю некоторого $A$. Рассмотрим последовательность средних 
\[
E(\mathbf{a})=\left(a_1, \frac{a_1+a_2}{2}, \frac{a_1+a_2+a_3}{3}, \ldots, \frac{a_1+\ldots+a_k}{k},\ldots\right)
\]
Фиксируем произвольный интервал $(\alpha,\beta)$ и подсчитаем, сколько раз последовательность средних пересекает этот интервал.

\begin{center}
\includegraphics[scale=1]{ergodic-1.pdf}
\end{center}

Можно проиллюстрировать пересечения картинкой. Построим ломаную из начала координат, каждый раз сдвигаясь на $1$ вправо и добавляя по вертикали очередное $a_i$. Тогда среднее значение $(a_1+\ldots+a_i)/i$ соответствует наклону луча, ведущего из начала координат к $i$-й вершине ломаной. Нас интересует выход этого наклона за пределы отрезка $[\alpha,\beta]$, то есть выход вершин ломаной из отмеченного на рисунке угла (стороны угла имеют наклоны $\alpha$ и $\beta$). В нашем конкретном примере мы видим одно пересечение щели <<сверху вниз>> (второй отрезок) и одно пересечение <<снизу вверх>> (два следующих звена). Дальше пересечений нет (наклон не опустился ниже $\alpha$).

Количество пересечений в обе стороны отличается не более чем на $1$, так что не так важно, что выбрать. Мы будем считать пересечения снизу вверх, то есть искать возрастающие последовательности $l_1 < r_1 < l_2 < r_2<\ldots$ индексов, сгрупированные в пары $(l_i,r_i)$: в $l_i$ среднее меньше $\alpha$, в $r_i$ больше $\beta$. Каждую пару (переход от <<меньше>> к <<больше>>) считаем за одно пересечение. Обозначим через $C_n(\mathbf{a})$ максимальное число таких пар среди первых $n$ членов последовательности средних $E(\mathbf{a})$ --- на нашей картинке, скажем $C_3(\mathbf{a})=0$, а $C_4(\mathbf{a})=\ldots=C_7(\mathbf{a})=1$. С ростом $n$ величина $C_n(\mathbf{a})$ (для данной последовательности $\mathbf{a}$) может только увеличиться за счёт появления новых пар. Её предел при $n\to\infty$ (конечный или бесконечный) мы обозначим $C(\mathbf{a})$. Здесь мы используем ту же букву $C$, что и раньше, поскольку ранее определённое $C(x)$ равно $C(\mathbf{x})$ для $\mathbf{x}=(f(x), f(Tx), f(TTx),\ldots$).

Следующая лемма оценивает среднее количество пересечений, подсчитанное \emph{для сдвигов произвольной ограниченной последовательности} $\mathbf{a}=(a_1,a_2,\ldots)$.  Другими словами, положим
\[
\mathbf{a}^{(i)}=(a_i,a_{i+1}, a_{i+2},\ldots)
\]
и найдём среднее число пересечений для $N$ таких сдвигов, то есть
\[
\frac{C_n(\mathbf{a}^{(1)})+C_n(\mathbf{a}^{(2)})+\ldots+C_n(\mathbf{a}^{(N)})}{N}.
\]

\begin{lemma}
При $N\ge n$ эта величина не превосходит $c\frac{A+|\alpha|+|\beta|}{\beta-\alpha}$ \textup(где $c$ --- некоторая абсолютная константа, а $A$ --- верхняя граница для всех $|a_i|$\textup).
\end{lemma}

Условие $N\ge n$ означает, что усредняется достаточно много чисел. В знаменателе дроби стоит $\beta-\alpha$, что не удивительно: чем уже щель, тем больше может быть пересечений. В числителе стоит $A$ (верхняя граница для всех $|a_i|$); чем $A$ больше, тем круче может идти кривая и тем больше может быть пересечений. Можно ещё заметить, что если умножить $A$, $\alpha$ и $\beta$ на одно и то же число, то оценка не изменится (что тоже не удивительно: умножение всех $a_i$,  а также $\alpha$ и $\beta$ на одно и то же число не меняет количество пересечений).

На картинке переход от $\mathbf{a}^{(i)}$ к $\mathbf{a}^{(i+k)}$ соответствует сдвигу начала отсчёта на $k$ единиц вправо.
\begin{center}
\includegraphics[scale=1]{ergodic-2.pdf}
\end{center}
Здесь, скажем, нарисованы два треугольника, пересечения которых нас интересуют: начальный (одно пересечение снизу вверх) и сдвинутый на два шага (ни одного пересечения).

В этой терминологии лемма говорит о том, что если мы посчитаем число пересечений снизу вверх для одинаковых треугольников, помещённых в разные места ломаной, то в среднем много не получится.

Используя эту лемму, легко закончить доказательство эргодической теоремы. Нас интересовало среднее значение $\mathbb{E}_x C(x)$, где $C(x)=C(\mathbf{x})$ --- число пересечений для последовательности $\mathbf{x}=(f(x), f(Tx), f(TTx),\ldots)$, и нужно было показать, что это среднее значение конечно (и потому множество точек, где $C(x)$ бесконечно, имеет меру нуль).  Для этого достаточно показать, что соответствующее среднее (по~$x$) значение $C_n(\mathbf{x})$ ограничено константой, не зависящей от $n$ (и потом сослаться на теорему о монотонной сходимости под знаком интеграла). Ключевое наблюдение: \emph{в силу инвариантности меры при преобразовании $T$ то же самое среднее будет и у $C_n(Tx)=C(\mathbf{x}^{(2)})$, и у $C_n(TTx)=C(\mathbf{x}^{(3)})$,\ldots} Поэтому можно дополнительно усреднить по первым $N$ итерациям преобразования $T$. А теперь, переставив порядок усреднения, можно заметить, что мы усредняем средние вдоль траектории из леммы (взяв $N$ не меньше $n$ --- при дополнительном усреднении $N$ можно взять любым). Поэтому среднее значение $C_n(x)$ не превосходит $c(F+|\alpha|+|\beta|)/(\beta-\alpha)$, где $F$ ---верхняя граница для модуля функции $f$ (она же --- граница для всех членов последовательности $\mathbf{x}$), а  $(\alpha,\beta)$ --- интересующая нас щель.

Что и завершает доказательство эргодической теоремы (точнее, вывод её из леммы, которую мы докажем в следующем разделе).

\section{Доказательство леммы}

Нам надо оценить сверху среднее число пересечений для нескольких подряд идущих сдвигов произвольной последовательности. Как и при анализе амортизированной сложности алгоритмов, тут полезно ввести некоторый \emph{потенциал} и показать, что большое число пересечений приведёт к большому увеличению потенциала, а в целом при последовательных сдвигах потенциал остаётся ограниченным. 

Функция потенциала (Вьюгин называл её \emph{кумулятивной суммой}) может показаться несколько искусственной, поэтому мы начнём с некоторых простых замечаний. 

Путеводители по горным прогулкам часто указывают не только перепад высот (разницу высот между концом и началом пути), но и общую высоту подъёмов на пути, то есть сумму перепадов высот для всех участков подъёма. Например, на картинке
\begin{center}
\includegraphics[scale=1]{ergodic-3.pdf}
\end{center}
отмечены три участка подъёма, и общий подъём равен сумме перепадов высот для каждого из трёх участков.

Суммарный подъём можно определить по-разному (скажем, можно взять минимальное представление нашей функции ограниченной вариации в виде разности двух монотонных, и взять изменение одной из них), но нам будет удобно использовать такое определение:
\begin{center}
\includegraphics[scale=1]{ergodic-4.pdf}
\end{center}
мы берём всевозможные последовательности подинтервалов 
\[
(l_1,r_1), (l_2,r_2),\ldots, (l_i,r_i)\qquad \text{с}\qquad 0\le l_1 < r_1 <  l_2 < r_2<\ldots<l_i<r_i\le n
\]
при любых $i$, для каждой из них вычисляем сумму
\[
h(r_1)-h(l_1)+h(r_2)-h(l_2)+\ldots+h(r_i)-h(l_i),
\]
и берём максимум этих сумм по всем последовательностям интервалов и  по всем $i$ (включая случай $i=0$, когда сумма считается нулевой). Это будет та же величина (потому что можно взять возрастающие участки, а невозрастающие можно вырезать, и станет только лучше).

Такой вариант определения подъёма легче приспособить для наших целей. Приспособление состоит в том, что мы вспоминаем о щели $(\alpha,\beta)$. Для правых концов $r_i$ вместо абсолютной высоты $h(r_i)$ мы берём превышение над верхней границей щели, то есть $h(r_i)-\beta r_i$, а для левых концов --- превышение над нижней границей щели, то есть $h(l_i)-\alpha l_i$. Другими словами, мы ищем максимум выражения
\[
[h(r_1)-\beta r_1] - [h(l_1)-\alpha l_1]+[h(r_2)-\beta r_2] - [h(l_2)-\alpha l_2]+\ldots+[h(r_i)-\beta r_i]-[h(l_i)-\alpha l_i],
\]
 В остальном определение не меняется (мы по прежнему берём все семейства непересекающихся подинтервалов рассматриваемого интервала индексов, и берём максимум).

Теперь нам надо, согласно объявленному плану, изучить два вопроса:
\begin{itemize}
\item как меняется эта величина (\emph{кумулятивная сумма}), когда мы сдвигаем начало отсчёта (другими словами, убираем первое звено в начале ломаной и добавляем новое звено к концу) и связать это изменение с числом пересечений щели;
\item почему кумулятивная сумма остаётся ограниченной.
\end{itemize}

Напомним обозначения: у нас есть последовательность $\mathbf{a}=(a_1,a_2,\ldots)$; ломаная начинается в точке $(0,0)$, затем идёт в точку $(1,a_1)$, потом в $(2,a_1+a_2)$, и так далее до $(n, a_1+a_2+\ldots+a_n)$. Нас интересует число пересечений щели $(\alpha,\beta)$ на отрезке $[0,n]$ (первые $n$ звеньев ломаной), а также кумулятивные суммы на отрезках $[0,n]$ и $[1,n+1]$. (Под кумулятивной суммой на $[1,n+1]$ мы понимаем кумулятивную сумму для последовательности $(a_2,a_3,\ldots)$ на отрезке $[0.n]$: другими словами, мы перемещаем начало отсчёта в конец первого отрезка ломаной, и добавляем справа ещё один отрезок с некоторым наклоном~$a_n$)

Что происходит при переходе от $[0,n]$ к $[1,n+1]$? Посмотрим на набор отрезков $0\le l_1 < r_1 <l_2 <r_2<\ldots<l_s<r_s\le n$, на котором достигается максимум (равный первой кумулятивной сумме). Предположим пока, что $l_1>0$. Тогда все эти отрезки годятся (в новой системе отсчёта, то есть с уменьшенными на $1$ концами) и для второй кумулятивной суммы, но соответствующие слагаемые немного изменятся, потому что вычитаемые в квадратных скобках уменьшатся на $\beta$ для правых концов и на $\alpha$ для левых концов. Каждый отрезок даст увеличение $(\beta-\alpha)$, так что кумулятивная сумма увеличится на $s(\beta-\alpha)$, а может, и больше (если появление последнего звена позволит выбрать лучшую точку максимума).

Особый случай возникает при $l_1=0$. Тогда придётся либо уменьшить первый отрезок (он будет $[0,r_1-1]$), либо вообще его удалить (если $r_1=1$, уменьшать некуда). Но эти изменения можно грубо оценить как $O(A+|\alpha|+|\beta|)$, где $A$ --- верхняя граница для всех $a_i$. В самом деле, вся область, где эти изменения происходят, имеет такой размер. Таким образом, гарантировано увеличение кумулятивной суммы на 
\[
s(\beta-\alpha)-O(A+|\alpha|+|\beta|),
\]
если максимум кумулятивной суммы достигается на наборе из $s$ отрезков.

Теперь покажем, что если пересечений щели $(\alpha, \beta)$ много, то число отрезков, при котором кумулятивная сумма достигает максимума (обозначенное нами через $s$) тоже должно быть большим. Мы докажем, что число пересечений щели не превосходит $O(s)$. При подсчёте пересечений мы чередуем точки, где наклон меньше $\alpha$, и точки, где наклон больше $\beta$. Пусть $u<v<w$ --- три соседние точки, и в первой и третьей наклон меньше $\alpha$, а во второй наклон больше $\beta$ (на рисунке ломаная изображена в виде кривой).

\begin{center}
\includegraphics[scale=1]{ergodic-5.pdf}
\end{center}

Сравним их с точками $l_1<r_1<l_2<r_2<\ldots< l_s<r_s$, дающими максимум в определении кумулятивной суммы. Эти последние разбивают прямую на интервалы, и мы утверждаем, что
\begin{center}
\emph{точки $u,v,w$ не могут попасть внутрь одного и того же интервала}.
\end{center}
В самом деле, если они попали в интервал $(l_i,r_i)$, то его можно заменить на  два интервала $[l_i,v]$ и $[w,r_i]$. От этого кумулятивная сумма увеличится, потому что к ней добавятся два положительных числа (соответствующие вертикальным отрезкам на рисунке --- показывающим, насколько $v$ больше верхней границы и насколько $w$ ниже нижней границы). Аналогичным образом, если все три точки разом попали в интервал $(r_i,l_{i+1})$, то можно добавить новый интервал $[u,v]$ между $[l_i,r_i]$ и $[l_{i+1},r_{i+1}]$, и по тем же причинам сумма увеличится. Можно провести аналогичное рассуждение и в том случае, когда первая и третья точки выше верхней границы, а вторая ниже нижней.

Отсюда следует, как мы и обещали, что для максимального набора отрезков (в кумулятивной сумме) их число $s$ не меньше числа пересечений с точностью до некоторого постоянного коэффициента (который легко указать явно,  но это нам не важно).

Второе обещанное утверждение: оценка кумулятивной суммы. Нижняя оценка: как мы уже говорили, она всегда неотрицательна. Для верхней оценки заметим, что если $\alpha\ge 0$, то кумулятивная сумма не больше максимальной длины участков подъёма (потому что в правых концах отрезков мы вычитаем больше, чем в левых: $\beta j \ge\alpha j \ge \alpha i$), и потому не превосходит $nA$. С другой стороны, кумулятивная сумма не изменится, если увеличить $\alpha$, $\beta$ и все члены последовательности $\mathbf{a}$ на одно и то же число, поэтому можно считать $\alpha$ неотрицательным, заменив $A$ на $A+|\alpha|$. Поэтому кумулятивная сумма на участке длины $n$ не превосходит $(A+|\alpha|+|\beta|)n$ --- мы добавили в скобку $\beta$ для симметрии.

Теперь можно собрать всё вместе, сложив оценки для изменения кумулятивной суммы для $N$ последовательных отрезков (для которых имеется $c_1,\ldots,c_N$ пересечений:
\[
\Omega(1)\sum_{i=1}^{N} c_i (\beta-\alpha)-NO(A+|\alpha|+|\beta|) \le O((A+|\alpha|+|\beta)n)
\]
Деля на $N$ и вспоминая, что $n\le N$ по предположению, заключаем, что среднее $N$ чисел пересечений $c_1,\ldots,c_N$ для $N$ последовательных участков ломаной не превосходит $O((A+|\alpha|+|\beta|)/(\beta-\alpha))$. Лемма доказана.

\section{Конструктивный вариант}

Как обычно, классические теоремы, верные почти всюду, имеют алгоритмические варианты, где <<почти всюду>> заменяется на <<для всех алгоритмически случайных точек>>, и обычно алгоритмическая случайность понимается в смысле определения Мартин-Лёфа.  Мы ограничимся пространством бесконечных последовательностей, где действует оператор сдвига (отбрасывания первого члена).  Мера $P$ на таком пространстве задаётся мерой интервалов (интервал $[x]$ --- это бесконечные последовательности, начинающиеся на слово $x$), то есть функцией $p(x)=P([x])$, определённой на двоичных словах и принимающей неотрицательные значения, для которой
\begin{itemize}
\item $p(\Lambda)=1$;
\item $p(x)=p(x0)+p(x1)$;
\end{itemize} 
Инвариантность меры относительно сдвигов (стационарность) означает, что в дополнение к этому для любого $x$ выполнено равенство
\begin{itemize}
\item $p(x)=p(0x)+p(1x)$.
\end{itemize}

Теперь надо наложить алгоритмические ограничения на меры и функции. Меру мы будем считать вычислимой (есть алгоритм, вычисляющий меру любого интервала с любой точностью). Функцию $f$ мы тоже будем считать вычислимой в аналогичном смысле (при этом аргумент даётся ей в качестве оракула). Из вычислимости следует непрерывность и ограниченность (пространство $X$ компактно).

\begin{theorem}[Вьюгин~\cite{vyugin}]
Пусть $P$~--- произвольная вычислимая стационарная мера. Пусть $f$~--- вычислимая функция на бесконечных последовательностях. Тогда для любой последовательности $\omega$, случайной по Мартин-Лёфу относительно~$P$, существует среднее функции $f$ вдоль траектории, начинающейся с $\omega$.
\end{theorem}

Другими словами, множество последовательностей, для которых вдоль траектории нет среднего, эффективно нулевое: по любому рациональному $\varepsilon>0$ можно эффективно перечислить семейство интервалов, содержащее все такие последовательности и имеющее $P$-меру меньше~$\varepsilon$.

В частности, если в качестве $f$ взять характеристическую функцию какого-то интервала, то видно, что всякая случайная последовательность имеет предел частот нулей и единиц,  пределы частот двухбитовых подслов и так далее.

Для доказательства достаточно заметить, что для всякой рациональной щели $(\alpha,\beta)$ и числа $m$ множество тех последовательностей, где на сдвигах средние имеют $m$  или больше пересечений, является эффективно открытым (напомним, что $f$ вычислима с любой заданной точностью, а в условии пересечения можно взять строгие неравенства <<меньше $\alpha$>> и <<больше $\beta$>>), и мера упомянутого множества не больше $O(1/m)$ (с константой, зависящей от максимума модуля функции $f$ и щели).

Это рассуждение доказывает более сильный результат, использующий понятие случайности относительно эффективно замкнутого класса мер (определение и обсуждение свойств см. в~\cite{survey}):

\begin{theorem}[Вьюгин, вариант с классом мер]
Пусть $f$~--- вычислимая функция на бесконечных последовательностях. Тогда для любой последовательности $\omega$, случайной по Мартин-Лёфу относительно класса стационарных мер, существует среднее вдоль траектории, начинающейся с $\omega$.
\end{theorem}

В терминах покрытий это означает, что по любому рациональному $\varepsilon>0$ можно эффективно перечислить семейство интервалов, содержащее все плохие последовательности и имеющее меру меньше $\varepsilon$ \emph{относительно любой стационарной меры $P$}. Разница в том, что мы строим (и построили) семейство интервалов, которое годится для всех мер. 

\section{Вопросы}

1. Можно ли доказать лемму более непосредственно, не вводя довольно странной кумулятивной суммы?

Есть даже совсем простая формулировка в частном случае. Пусть $X$~--- двоичное слово. Будем следить за частотой единиц в начальных отрезках, и смотреть, как она колеблется и сколько раз становится то меньше $0.4$, то больше $0.6$ (границы выбраны для примера).  Это число назовём <<числом колебаний>>. Утверждение: \emph{каково бы ни было $n$, если мы возьмём в достаточно длинном слове все подслова длины $n$, то среднее число колебаний в них не превосходит $1000$} (константа выбрана наугад, но, наверно, должно хватить). Заметим, что при больш\'{и}х $n$ число колебаний в одном подслове длины $n$ может быть очень большим (но тогда в сдвигах должно быть меньше).

2. Утверждение теоремы Биркгофа верно для функций с векторными значениями, скажем, в $\mathbb{R}^2$, и формально следует из скалярного варианта. Однако в двумерном случае идея подсчёта пересечений выглядит более странно, так как там нет естественных границ, переход через которые можно было бы подсчитывать. Может быть, можно как-то более инвариантно изложить доказательство в этом случае?

3. Эффективная эргодическая теорема для случайных по Мартин-Лёфу последовательностей верна и для перечислимых снизу $f$ (или попаданий в эффективно открытое множество $A$), \emph{если преобразование является эргодическим}, см.~\cite{hoyrup}. А можно ли обойтись без этого предположения при таком обобщении?

4.~Эргодическая теорема для аменабельных групп в частном случае плоскости гарантирует, что если дана мера $P$ на пространстве плоских конфигураций из нулей и единиц, которая инвариантна относительно горизонтальных и вертикальных сдвигов, то множество тех конфигураций, в которых предельная частота единиц (по возрастающим последовательностям квадратов с фиксированным центром) не существует, имеет $P$-меру нуль. Спрашивается, можно ли здесь сделать аналогичную штуку и построить открытое множество, содержащее все такие (не имеющие предела частот) последовательности, и имеющее малую $P$-меру \emph{для любой инвариантной меры $P$}. Кажется, что рассуждение с колебаниями тут никак не применимо. Связанный вопрос: верен ли эффективный вариант теоремы Биркгофа (как у Вьюгина) для двумерных случайных последовательностей?